\documentclass[12pt]{amsart}

\def\@typesizes{%
       \or{5}{6.5}\or{6}{7.5}\or{7}{8.5}\or{8}{11}\or{9}{12}%
       \or{10}{13}% normalsize
       \or{\@xipt}{14}\or{\@xiipt}{15}\or{\@xivpt}{18}%
       \or{\@xviipt}{20}\or{\@xxpt}{24}}

\usepackage{amsthm}
\usepackage{amsmath}
\usepackage{amssymb}
\usepackage{graphicx}
\usepackage{enumerate}
\allowdisplaybreaks
\numberwithin{equation}{section}
\numberwithin{figure}{section}

\theoremstyle{plain}
\newtheorem{theorem}{ Theorem}[section]
\newtheorem{proposition}[theorem]{ Proposition}
\newtheorem{lemma}[theorem]{ Lemma}
\newtheorem{corollary}[theorem]{ Corollary}
\newtheorem{example}[theorem]{ Example}
\newtheorem{remark}[theorem]{ Remark}
\newtheorem{definition}[theorem]{ Definition}
\newtheorem{conjecture}{ Conjecture}

\def\BET{\begin{theorem}}
\def\ENT{\end{theorem}}
\def\BEP{\begin{proposition}}
\def\ENP{\end{proposition}}
\def\BEL{\begin{lemma}}
\def\ENL{\end{lemma}}
\def\BEC{\begin{corollary}}
\def\ENC{\end{corollary}}
\def\BEE{\begin{example} \rm}
\def\ENE{\end{example}}
\def\BER{\begin{remark} \rm}
\def\ENR{\end{remark}}
\def\BED{\begin{definition} \rm}
\def\END{\end{definition}}
\def\BECJ{\begin{conjecture}}
\def\ENCJ{\end{conjecture}}

\def\bea{\begin{eqnarray}}
\def\eea{\end{eqnarray}}

\def\beas{\begin{eqnarray*}}
\def\eeas{\end{eqnarray*}}

\def\beq{\begin{equation}}
\def\eeq{\end{equation}}

\def\beal{\begin{align*}}

\def\eeal{ \end{align*} }

\def\roweq{\nonumber \\ &=& }

\def\rowpl{\nonumber \\ & & +  }
\def\rowmi{\nonumber \\ & &  -  }

\def\bbC{{\mathbb C}}
\def\bbD{{\mathbb D}}

\def\bbN{{\mathbb N}}

\def\bbR{{\mathbb R}}

\begin{document}

\title[Solid hulls of weighted Banach spaces of analytic functions]{Solid hulls of weighted Banach spaces of analytic functions on the unit disc with exponential weights  }

\author{Jos\'{e} Bonet  and Jari Taskinen }

\begin{abstract}
Let $v(r)=\exp(-a/(1-r)^b)$ with $a>0$ and $0 < b \leq 2$ be an exponential weight on the unit disc. We study the solid hull of its associated weighted Banach space $H_v^\infty(\bbD)$ of all the analytic functions $f$ on the unit disc such that $v|f|$ is bounded.
\end{abstract}

\maketitle
\vspace{.4cm}

\section{Introduction}
\label{sec1}

Recently, the authors characterized in \cite{BT} the solid hulls of weighted $H^\infty$-type
Banach spaces of  entire functions for a large class of weight functions $v$.
An analytic function $f(z)=\sum_{n=0}^\infty a_n z^n$ on the disc is identified with the sequence of its Taylor coefficients $(a_n)_{n=0}^\infty$. Let $A$ be vector spaces of complex sequences containing the space of all the sequences with finitely many non-zero coordinates. The space $A$ is \textit{solid} if $a=(a_n) \in A$ and $|b_n| \leq |a_n|$ for each $n$ implies $b=(b_n) \in A$. The \textit{solid hull of $A$} is
$$S(A):= \{ (c_n) \, : \, \exists (a_n) \in A \ \mbox{such that} \ |c_n| \leq |a_n| \ \forall n \in \bbN \}. $$
It coincides with the smallest solid space containing $A$. See \cite{AS}.

The aim of this paper is to extend the results of \cite{BT} for
the corresponding spaces on the open unit disc $\bbD$. Accordingly, we study
Banach spaces  $H_v^\infty(\bbD)$ of analytic functions $f:\bbD \to \bbC$
such that $||f||_v:= \sup_{z \in \bbD} v(z)|f(z)| < \infty$. We show in Theorem \ref{th4.5} that the solid hull  for  $v(z) =
\exp( -1/(1-|z|)), z \in
\bbD,$ consists of complex sequences $(b_m)_{m=0}^\infty$ such that
\begin{equation}
\sup_{n \in \mathbb{N}}
\sum_{m=n^4 + 1 }^{(n+1)^4} |b_m|^2
\exp( -2n^2 ) \Big(1  - \frac1{n^2} \Big)^{2m}  %%\Big)^{1/2}
< \infty . \nonumber
\end{equation}
We also formulate a general Theorem \ref{th2.1}, which contains the characterization of
the solid hulls for a large  class of  weights.
This class of weights includes those satisfying condition (B) of \cite{L}.
Theorem \ref{th4.5} is used to determine space of multipliers from $H_v^\infty(\bbD)$ into $\ell^p, 1 \leq p \leq \infty,$ in Proposition \ref{propmult}.

Surprisingly enough, we will encounter a technical difficulty, which makes the calculation
of the solid hulls for weights $v(z) = \exp(-a/(1-|z|)^b)$ on the disc  more complicated than
those for the somewhat analogous weights $v(z) = \exp(-a|z|^b)$ in the plane. The latter were
successfully treated in \cite{BT} for all $b > 0$ and $a > 0$, but
it seems to the authors that in the case of the disc the calculation of certain numerical
sequences requires approximate solutions of some numerical equations, which can only be
done relatively simply for small enough $b$: for the above mentioned weights, we only
complete the calculation in the case $0 < b \leq 2$.

Bennet, Stegenga and Timoney in their paper \cite{BST} determined the solid hull and the solid core of the weighted spaces $H^\infty_v(\bbD)$  in the case the weight $v$ is doubling. Exponential weights $v(r)=\exp(-a/(1-r)^b)$ with $a,b>0$ are not doubling.
Not much seems to be known about multipliers and solid hulls of weighted spaces of analytic functions on the unit disc in the case of exponential weights. Hadamard multipliers of certain weighted space $H^1_a(\alpha), \alpha >0,$ were completely described by Dostani\'c in \cite{Dos} (see also Chapter 13 in \cite{JVA}). Other aspects of weighted spaces of analytic functions on the unit disc with exponential weights, like integration operators or Bergman projections, have been investigated recently by Constantin, Dostanic, Pau, Pavlovi\'c, Pel\'aez and R\"atty\"a, among others; see \cite{CP}, \cite{Dos2}, \cite{PP}, \cite{Pav} and \cite{PR}. The solid hull and multipliers on spaces of analytic functions on the disc has been investigated by many authors. In addition to  \cite{BST}, we mention for example  \cite{AS}, \cite{JP}, the books \cite{JVA} and \cite{Pav-book} and the many references therein.

Spaces of type $H_v^\infty(\bbC)$ and $H^\infty_v(\bbD)$ appear
in the study of growth conditions of
analytic functions and have been investigated in various articles since
the work of Shields and Williams, see {\it e.g.} \cite{BBG},\cite{BBT},
\cite{L}, \cite{SW} and the references therein.

A {\it weight} $v$   is a continuous function  $v: [0, 1[ \to ]0,  \infty [$, which is non-increasing on $[0,1[$ and satisfies
$\lim_{r \rightarrow 1} v(r)=0$. We extend $v$ to $\bbD$ by $v(z):= v(|z|)$. For
such a weight, the {\it weighted Banach space of analytic functions}
is denoted by $H^\infty_v(\bbD)$ and its norm by  $\Vert \cdot  \Vert_v .$  For an analytic
function $f \in H(\bbD)$, we denote $M(f,r):= \max\{|f(z)| \ | \ |z|=r\}$. Using the notation $O$ and $o$ of Landau, $f \in H_v^\infty(\bbD)$ if and only if $M(f,r)=O(1/v(r)), r \rightarrow 1$.

\section{The results.}
\label{sec2}
As mentioned in the introduction, the solid hull of the weighted Banach space $H^\infty_v(\bbD)$ of holomorphic functions on the unit disc $\bbD$ when the weight $v$ is doubling, i.e.
$$
\sup_{0<s\leq 1} \frac{v(1-s)}{v(1- 2^{-1}s)} < \infty,
$$
was determined in \cite{BST}.
The doubling condition is equivalent to the condition $(*)$ of \cite{L}, see
Example 2.4 of the citation. This condition appears also e.g.  in \cite[Theorem 3.2]{BDLT}.
Examples of weights on the disc which are not doubling are given by
$v(r):= \exp(-a/(1-r)^b)$, $a>0$, $b>0$.
In this section we investigate the solid hull of $H^\infty_v(\bbD)$ for
a large class of non doubling weights.

Given a continuous, radial weight $v$ and $m>0$, we denote by $r_m$ the global maximum point
in $[0,1[$  of the function $r \mapsto r^m v(r)$. For $0<m<M$ we define
\begin{equation}
A(m,M):= \left(\frac{r_{m}}{r_{M}}\right)^{m} \frac{v(r_{m})}{v(r_{M})} \ \
\ \mbox{and}  \ \ \
B(m,M):= \left(\frac{r_{M}}{r_{m}}\right)^{M} \frac{v(r_{M})}{v(r_{m})}.
\label{28}
\end{equation}
A word by word repetition of the proof of Theorem 2.5. of
\cite{BT} yields the next general result; notice that the argument uses the results of
\cite{L}, which are also established for holomorphic functions on the unit disc.

\begin{theorem}\label{th2.1}
Let $v$ be a radial weight on $\mathbb{D}$. Let $0<m_1 <m_2<...$ be a sequence with
$\lim_{n \rightarrow \infty} m_n = \infty$, such that for some constants
$2 < k \leq K$ we have
\begin{equation}
k \leq A(n):= A(m_n,m_{n+1}) \leq K \ \ , \ \
k \leq B(n):=B(m_n,m_{n+1}) \leq K
\label{225}
\end{equation}
for each $n \in \bbN$. Then, the solid hull of $H^\infty_v(\bbD)$ is
\begin{equation}
S(H^\infty_v(\bbD)):= \Big\{ (b_m)_{m=0}^\infty \, : \,
\sup_{n} v(r_{ m_{n}})
\Big( \sum_{m=m_n+1}^{m_{n+1}} |b_m|^2
 r_{ m_n}^{2m} \Big)^{1/2} < \infty \  \Big\} . \label{227}
\end{equation}
\end{theorem}

Proceeding as in \cite[Remark 2.7]{BT} one can show that every weight on the disc satisfying condition (B) of Lusky  \cite{L} satisfies the assumptions of Theorem \ref{th2.1}.

The main result of this paper is the calculation of the solid hull for some of the
most  usual non-doubling weights on the unit disc, namely the weights
\bea
v(r)= \exp(- a/(1-r)^b) ,   \label{227k}
\eea
where $a>0$, $0 < b \leq 2$ are constants.

We write
\bea
\alpha = 2 + \frac2b, \ \ \beta = \frac1{1+b}, \ \  G=(ab)^\beta .
\label{228}
\eea

\begin{theorem}\label{th4.5}
If $0 < b \leq 1$, then the solid hull $S(H^\infty_v(\bbD))$ of $H^\infty_v(\bbD)$ consists
of sequences $(b_m)_{m=0}^\infty$ satisfying
\begin{equation}
\sup_{n \in \bbN} \exp\big(  - 2 aG^{-\beta} n^2   \big)
\!\!\!\!\!
\sum_{m=n^{2+2/b} +1 }^{(n+1)^{2+2/b}}
 \!\!\!\!\! |b_m|^2
\Big( 1 - \frac{G}{n^{2/b} } \Big)^{2m} < \infty .
 \label{230}
\end{equation}
In the case $1 < b \leq  2$, the condition for the solid hull is
 \begin{equation}
\sup_{n \in \bbN} \exp\big( -  2 aG^{-\beta} n^2 + 2 \beta (ab)^{2\beta} n  \big)  \!\!\!\!\!
\sum_{m=n^{2+2/b} +1 }^{(n+1)^{2+2/b}}  \!\!\!\!\! |b_m|^2
\Big(1 - \frac{G}{n^{2/b}} + \frac{\beta G^2 }{n^{4/b}} \Big)^{2m} < \infty .
 \label{231}
\end{equation}
\end{theorem}

In particular, for $v(r) = \exp(-1/(1-r))$ the solid hull is
\begin{equation}
\Big\{ (b_m)_{m=0}^\infty \, : \,
\sup_{n} \exp ( -2n^2)
\sum_{m=n^4 +1 }^{(n+1)^4} |b_m|^2 \Big( 1 - \frac1{n^2} \Big)^{2m}
< \infty \  \Big\} ,  \label{232}
\end{equation}
and for $v(r) = \exp(-1/(1-r)^2 )$ it  is (!)
\begin{equation}
\Big\{ (b_m)_{m=0}^\infty \, : \,
\sup_{n \in \bbN} \exp \Big( -2^{8/9}    n^2 + \frac{2^{5/3}}3 n \Big)
\sum_{m=n^3 +1 }^{(n+1)^3} |b_m|^2 \Big( 1 - \frac1{2^{1/3} n} +  \frac1{3 \cdot  2^{2/3} n^2} \Big)^{2m}
< \infty \  \Big\} .  \label{232p}
\end{equation}

The proof requires the choice  of the numbers $m_n$ such that
the condition \eqref{225} is satisfied. This  turns out to
be quite technical and will be done in the next section. The main complication is
that, contrary to the case of weights $\exp(-a r^p)$ of the paper \cite{BT},
the maximum points $r_m$ cannot be solved explicitly (see \eqref{235t}),
and one has to treat only
approximations of them (see Lemma \ref{lem4.2}) in the calculations.

The proof of
Theorem \ref{th4.5} will be completed  in Section \ref{sec4}, where we also
make some comments on the case $b >2$.
Consequences of Theorem \ref{th4.5} for multipliers from $H^{\infty}_v(\bbD)$ into
$\ell_p, 1 \leq p \leq \infty$ are given in Section \ref{secmult}.

\section{Calculation of some numerical sequences.}
\label{sec3}

For this section, let $v$ be as in \eqref{227k}.
The next lemma will be crucial for the proof of
Theorem \ref{th4.5}, and its proof will occupy
the whole  section.

\begin{lemma}
\label{lem4}
For the weight $v$, the quantities
$A (n)$ and $B(n)$ satisfy \eqref{225}, if $m_n$ is chosen to be
\begin{equation}
m_n= S n^{2 + \frac{2}{b}} =: S n^\alpha ,  \label{233}
\end{equation}
where  $S = (16a)^{-1} $ for $b =1$,  $S =  ba^{-1/b}\alpha^{-1 - 1/b} $ for $b <2$,
and $S= \max \{ 2a, 2a^{-1/2} \}$  for $b=2$.

\end{lemma}
We will need  several times the Taylor expansions
\bea
(1+ x)^{c} & = & 1 +c x + \frac12 c(c-1)x^2 + O(x^3) \label{233a} \\
\log (1 + x ) &=& x - \frac12 x^2  + \frac13 x^3 + O(x^4), \label{233b}
\eea
valid for $|x| < 1$ and $c \in \bbR$.
We start by the
following estimate. Recall that $\alpha, \beta$ and $G$ are defined in (\ref{228}).

\begin{lemma}
\label{lem4.2}
Given $m \geq \max \{ 1, ab\}$, the global maximum point $r_m$ of the function $r \mapsto
r^mv(r)$ satisfies the estimate
\begin{equation}
1  - \frac{G}{ m^{\beta}}  + \frac{ \beta G^{2}}{m^{2 \beta} }
- \frac{C}{ m^{3 \beta}}
\leq r_m
\leq 1  - \frac{ G}{ m^{\beta}}
+ \frac{\beta G^2}{m^{2 \beta} }
+ \frac{C}{ m^{3 \beta}} .
\label{235}
\end{equation}
\end{lemma}

Proof.
A simple calculation shows that $r_m$ satisfies the equation
\begin{equation}
m (1-r_m)^{b+1} - ab r_m = 0 .  \label{235t}
\end{equation}
It is obvious that for some fixed $0 < \epsilon \leq 1/2$, the equation \eqref{235t}
has a solution $r_m \in [\epsilon,1[$, for all $m \geq 1 $. Writing $\delta = 1 -r_m \in
]0,1]$, \eqref{235t} is equivalent with
\begin{equation}
m^{\frac{1}{b+1}} \delta = (ab)^{\frac{1}{b+1}}  (1 - \delta)^{\frac{1}{b+1}}
= (ab)^{\frac{1}{b+1}} \Big( 1 - \frac{1}{b+1} \delta -
 \frac{b}{2(b+ 1)^2}  \delta^2 - +  \ldots \Big),
\nonumber
\end{equation}
or, by \eqref{228},
\begin{equation}
\delta =  \frac{G}{ m^{\beta}} \big( 1 - \beta \delta  %%- b \beta^2 \delta^2 +
+ \delta^2 F(a,b, \delta) \big)
\label{238}
\end{equation}
where the expression $F(a,b,\delta)$ is for fixed $a$, $b$, uniformly
bounded in $\delta \in ]0,1-\epsilon]$. (The choice of $\epsilon$ is needed here
as $F$ would not be uniformly bounded for  $\delta \in ]0,1]$.)
We already noticed  that \eqref{238} has a solution
$\delta \in [0,1-\epsilon]$. Setting the estimate $\delta \leq 1$ to the right hand side
of \eqref{238} yields the bound $\delta \leq Cm^{-\beta}$. Since $m \geq ab$ by
assumption, we have $Gm^{-\beta} \leq 1$ and thus the  solution $x$
(\eqref{239u}, below) of the equation
\begin{equation}
x =  \frac{G}{ m^{\beta}} ( 1 - \beta x)  \label{239} %%  + b \beta x^2 \big)
\end{equation}
satisfies $0 < x < 1$;   putting the bound $x \leq 1$ to the right hand side of
\eqref{239} yields   $x \leq Cm^{-\beta}$.  Then, subtracting \eqref{239} from
\eqref{238} and using the triangle inequality,
\begin{equation}
|\delta - x | \leq  \frac{G \beta}{ m^{\beta}} |\delta - x | + \frac{G \delta^2}{ m^{\beta}} |F| .
\label{239a}
\end{equation}
But we have  $|x - \delta | \leq x + \delta \leq C m^{-  \beta} $, and putting this
to the right hand side of \eqref{239a} implies $|\delta - x | \leq Cm^{-2\beta}$.
Substituting this once more in \eqref{239a} yields $|\delta - x | \leq Cm^{-3\beta}$.
This implies \eqref{235}, since the solution of \eqref{239} is
\begin{equation}
x= \frac{Gm^{-\beta}}{1+ \beta Gm^{-\beta}} = Gm^{-\beta} - \beta G^2 m^{-2\beta}
+ O(m^{-3\beta}) .  \ \ \Box \label{239u}
\end{equation}

From now on we assume that $m_n$  is defined for all  $n \in \bbN$
as in \eqref{233}. By \eqref{235} we can write for all $n$ such that $m_n \geq
\max \{ 1, ab\}$,
\bea
1 -r_{m_n} = Gm_n^{-\beta} - \beta G^2 m_n^{-2 \beta} + R(m_n)  \label{240-}
\eea
where $|R(m_n) | \leq Cm_n^{-3 \beta}$.

\begin{lemma}
\label{lem4.3}
Let  $0 < b \leq 2$. We have for all large enough $n \in \bbN$
\begin{equation}
\log \left( \frac{v(r_{m_{n+1}})}{v (r_{m_n})}
\right)  =
\left\{
\begin{array}{ll}
- {\displaystyle \frac{2G}b} S^{b\beta} n
-  {\displaystyle \frac{G}{b}} S^{b \beta}  + O(n^{-1} + n^{1-2/b}) ,
 \ &  \mbox{if} \ b < 2  \\
\ & \ \\
- G S^{2/3} n
-  {\displaystyle \frac{G}2} S^{2/3} +  {\displaystyle \frac13} G^2 S^{1/3} + O(n^{-1}) ,  &  \mbox{if} \ b=2.
\end{array}
\right.
\label{240}
\end{equation}
\end{lemma}
Notice that in the case  $b=2$ we have $b \beta=2/3$ and more importantly,
the constant term (of order $n^0$) is not the same as what would be gotten from
the formula in the case $b < 2$.

\bigskip

Proof. Using \eqref{240-}  and \eqref{233a}, and assuming that
$n$ is so large that $Gm_n^{-\beta} < 1$,  we get
\bea
v(r_{m_n})  &= &\exp\Big(-a \Big( Gm_n^{-\beta} - \beta G^2 m_n^{-2 \beta} +
R(m_n) \Big)^{-b} \Big)
\roweq
\exp\Big(-a G^{-b}m_n^{b\beta} \Big( 1 - \beta G m_n^{- \beta} +
G^{-1} m_n^{ \beta}R(m_n) \Big)^{-b} \Big)
\roweq
\exp\big(-a G^{-b}m_n^{b\beta}  + ab \beta G^{1-b} m_n^{(b-1) \beta } +
\tilde R_n  \big)   \label{241x}
\eea
where $|\tilde R_n| \leq C m_n^{\beta (b-2) }$.
We thus get, taking into account that $m_n= Sn^\alpha$ and using
\eqref{233a} again,
\beas
& & \log \left( \frac{v(r_{m_{n+1}})}{v (r_{m_n})}
\right)
\roweq
-a G^{-b}\big(m_{n+1}^{b \beta} -m_n^{b\beta}\big)  + ab \beta G^{1-b}
\big(m_{n+1}^{(b-1)\beta} -  m_n^{(b-1) \beta }\big) +
\tilde R_{n+1} - \tilde R_n
\roweq
-a G^{-b} S^{b \beta} \big( (n+1)^{b \alpha \beta} -n^{b\alpha \beta}\big)  +
\rowpl
ab \beta G^{1-b}
 S^{(b-1) \beta} \big((n+1)^{(b-1)\alpha \beta} -  n^{(b-1)\alpha \beta }\big)
+ \tilde R_{n+1} - \tilde R_n
\roweq
-aG^{-b}  S^{b \beta} b \alpha \beta n^{b \alpha  \beta-1}
-aG^{-b}S^{b \beta} \frac12  b\alpha \beta (b \alpha  \beta -1)
 n^{b \alpha \beta-2} + O( n^{ b\alpha \beta-3})
\rowpl
 ab \beta G^{1-b}  S^{(b-1) \beta}
 (b-1)\alpha \beta n^{(b-1) \alpha\beta - 1} +O(n^{(b-1)\alpha\beta - 2}).
%\label{240a}
\eeas
This is simplified to the claimed form by observing that
\bea
& & ab G^{-b} = (ab)^{1-\frac{b}{1+b}} = (ab)^\beta = G  \label{257} \\
& &  \alpha \beta = 2(1+b)b^{-1}(b+1)^{-1}  = 2/b \ \ , \ \
  b \alpha \beta = 2  \label{258}
\eea
and moreover
\begin{equation}
 (b-1) \alpha \beta -1 =  b \alpha \beta - 1 - \alpha \beta
= 1 - 2b^{-1}
\left\{
\begin{array}{ll}
< 0 , & b < 2 , \nonumber \\
= 0 , &b=2 .
\end{array}
 \nonumber
\right.
\end{equation}

\begin{lemma}
\label{lem4.4}
If $0 < b < 2$, we have for all large enough $n \in \mathbb{N}$,
\bea
& & m_{n+1} \log  \Big( \frac{r_{m_{n+1}}}{r_{m_n}} \Big)
\roweq
\frac{2G}b S^{b\beta} n
+  \frac{ G}{b}  \Big(3+ \frac{2}{b }\Big)S^{b \beta}
+ O(n^{-1} + n^{ 1-2/b }) .
\label{940}
\eea
If $b=2$,  then
\bea
& & m_{n+1} \log  \Big( \frac{r_{m_{n+1}}}{r_{m_n}} \Big)
= G  S^{2/3 } n
+  2 GS^{2/3 }  + \frac13  G^2 S^{1/3} + O(n^{-1} ) .
\label{940b}
\eea
\end{lemma}

Proof. We have again by \eqref{240-}, \eqref{233b}, for large enough $n$,
\beas
& & m_{n+1} \log  \Big( \frac{r_{m_{n+1}}}{r_{m_n}} \Big)
\roweq
m_{n+1}\bigg( \log \Big( 1 - Gm_{n+1}^{-\beta}  +\beta
G^2  m_{n+1}^{-2 \beta}  +  R(m_{n+1})  \Big)
\rowmi
\log \Big( 1 - Gm_{n}^{-\beta}  +\beta  G^2  m_{n}^{-2 \beta}
+   R(m_n)\Big) \bigg)
\roweq
m_{n+1} \Big( -  G(m_{n+1}^{-\beta} -m_{n}^{-\beta})
+ \big( \beta - \frac12 \big)  G^2 ( m_{n+1}^{-2 \beta} - m_{n}^{-2 \beta} )
+ O(m_n^{-3 \beta})  \Big).
%%\label{940a}
\eeas
Here, keeping in mind that  $m_n = Sn ^\alpha$, the term
\bea
m_{n+1} \big( \beta - \frac12\big)  G^2 \big(  m_{n+1}^{-2 \beta} - m_{n}^{-2 \beta} \big)
\label{943}
\eea
is of degree  $\alpha - 2 \beta \alpha -1 =  1 - 2/b $ with respect to $n$,
and this number is negative, if and only if $b <2$. So, in the case $b < 2$
we obtain using $1 - \beta = b \beta$ and \eqref{233a}
\bea
& &
- m_{n+1} G(m_{n+1}^{-\beta} -m_{n}^{-\beta})
\roweq
- S(n+1)^\alpha G S^{- \beta}\big( - \alpha \beta n^{-\alpha \beta -1}
+ \frac12 \alpha \beta (\alpha \beta +1) n^{-\alpha \beta -2} + O(n^{-\alpha \beta -3}
\big)
\roweq
-  G  S^{b\beta}\big(n^{\alpha} + \alpha n^{\alpha -1} + O(n^{\alpha -2}) \big)
\nonumber \\
& & \times \ \big( -\alpha \beta n^{-\alpha \beta -1}
+ \frac12 \alpha \beta (\alpha \beta +1) n^{-\alpha \beta -2}
+O(n^{- \alpha \beta - 3}) \big) . \label{947}
\eea
As for the exponents, notice that
\bea
\alpha - \alpha \beta -1 =
2+ 2b^{-1} -2b^{-1} -1 = 1 .
\eea
The coefficient of $n$, respectively, $n^0$,  thus equals
\bea
G  S^{b\beta} \alpha \beta ,  \ \ \mbox{resp.} \ \
G  S^{b\beta}\Big(- \frac12 \alpha \beta (\alpha \beta +1) +
\alpha^2 \beta\Big)  \label{948}
\eea
This yields the claim of the lemma for $b < 2$  by  using \eqref{258}.

In the case $b=2$ the term \eqref{943} equals
\beas
& & S(n+1)^\alpha \big( \beta - \frac12 \big)  G^2  S^{- 2\beta}\big( - 2 \alpha \beta
n^{- 2 \alpha \beta -1 }
+ O(n^{- 2 \alpha \beta - 2} ) \big)
\roweq
 \frac13 S^{1/3} G^2 + O(n^{-1}) .  %\label{949}
\eeas
Adding this to the previous case yields \eqref{940b}. \ \ $\Box$

\begin{lemma}
\label{lem4.5}
If $0 < b \leq 2$, we have for all large enough $n \in \mathbb{N}$,
\begin{equation}
m_{n} \log  \Big( \frac{r_{m_{n+1}}}{r_{m_n}} \Big)
= \left\{
\begin{array}{ll}
{\displaystyle \frac{2G}b S^{b\beta} n
 -  GS^{b \beta} \frac1b \Big( \frac2b +1\Big)  + O( n^{1-2/b})} , & \ b < 2 \\
\ & \ \\
{\displaystyle  G S^{2/3} n
 -   GS^{2/3} +  \frac13 S^{1/3} G^2 + O(n^{-1}) } , & \ b = 2
\end{array}
\right.
\label{950}
\end{equation}
\end{lemma}

Proof. One makes the obvious change  $m_{n+1} \to m_n$, or, $(n+1)^\alpha \to n^\alpha$
in the proof of Lemma \ref{lem4.4} and collects the coefficients
of the remaining terms in the same way as in the argument \eqref{947}
(one obtains \eqref{948} except for the term $\alpha^2 \beta$).
The case $b=2$ is proven in the same way: in addition to the omission of the
$\alpha^2 \beta$-term there are no other  changes.
\ \ $\Box$

\bigskip

Proof of Lemma \ref{lem4}. Let first $b <2$.
We consider the quantity
$$
\log B(n) = m_{n+1} \log \Big( \frac{r_{m_{n+1}}}{r_{m_n}}  \Big)
+ \log  \Big( \frac{v(r_{m_{n+1}})}{v(r_{m_n})}  \Big)
$$
and first observe that the coefficients of the
term with $n$ are the opposite numbers in \eqref{240} and \eqref{940}.
The sum of  the coefficients of the term $n^0$ in \eqref{240}
and \eqref{940} is
\beas
 GS^{b \beta} \Big( - \frac1b + \frac3b + \frac{2}{b^2}\Big)
\eeas
so that we get
\bea
 \log B (n) & =&  GS^{b \beta} \Big( \frac{2}{b} + \frac{2}{b^2}\Big) +
O(n^{-1} + n^{1-2/b})
\roweq
 GS^{b \beta}  \frac{\alpha }{b}  + O(n^{-1} + n^{1-2/b})  . \label{980}
\eea
The required property \eqref{225} follows for $B(n)$ by choosing $S$ large enough
so that the constant term on the right of \eqref{980}
is at least 1.
We choose $S$ such that
\beas
GS^{b \beta} = \frac{b}{\alpha } \Rightarrow
S =  ba^{-1/b}\alpha^{-1 - 1/b} .  % \label{981}
\eeas
If $b =1$, we have  $\alpha =4$ so that
\beas
S= \frac1{16 a} . % \label{981a}
\eeas

The same calculation, using Lemma \ref{lem4.5} instead of
Lemma \ref{lem4.4}, yields (notice the order of the  numerator and denominator
in $A(n)$)
\bea
\log A(n) &=&
GS^{b \beta} \Big(  \frac1b + \frac1b \Big( \frac2b +1\Big) \Big) +
O(n^{-1} + n^{1-2/b})
\roweq
GS^{b \beta}\frac{\alpha}{b}  +
O(n^{-1} + n^{1-2/b})   ,  \label{982}
\eea
and we get the desired conclusion for $A(n)$ by the same
choice of $S$ as above.

Finally, if $b=2$,  we have instead of \eqref{982}
\bea
\log A(n) &=&
GS^{1/3} \Big(\frac32 S^{1/3} - \frac{2}{3}G \Big) +
O(n^{-1})  ,  \label{986}
\eea
where $G = (2a)^{1/3}$. Choosing
\bea
S = \max \{ 2a , 2 a^{-1/2} \}   \label{985}
\eea
the leading term in \eqref{986} has the estimate
\beas
& & GS^{2/3} \Big( \frac{3}{2} - \frac23 \frac{ G }{ S^{1/3}} \Big)
\geq
  GS^{2/3} \Big( \frac{3}{2} - \Big( \frac{2a}{S}\Big)^{1/3}\Big)
\geq
 \frac{1}{2} (2aS^2)^{1/3}  \geq  1 ,
\eeas
hence, \eqref{225} follows for $A(n)$, if $n$ is large enough.

We have  instead of \eqref{980} the estimate (cf. \eqref{240} and  \eqref{940b})
\beas
\ \ \log B (n)=
GS^{1/3} \Big( \frac{3}{2}S^{1/3} + \frac23  G \Big) +
O(n^{-1}) ,
%\label{984}
\eeas
and we thus see that the choice \eqref{985} is enough to guarantee
that $\log B(n) \geq 1 + O(n^{-1})$.   \ \ $\Box$

\bigskip

\section{Proof of Theorem \ref{th4.5}.}
\label{sec4}

We choose $m_n$ according to \eqref{233} for all $n\in\bbN$.
Theorem \ref{th4.5} follows in principle from
Theorem \ref{th2.1} and Lemma \ref{lem4}, but we need to be careful
to use accurate enough approximations of $r_{m_n}$.
If $b \leq 1$, we observe that in \eqref{241x}, the exponent of $n^{(b-1)\alpha \beta} = n^{2-2/b}$ in
$m_n^{(b-1)\beta}$ is at most 0 and $|\tilde R_n|$ is
clearly bounded by a constant, so we get for all $n \in \bbN$ using
 \eqref{258},
\bea
c_1 \exp(-aG^{-b} n^2) \leq v(r_{m_n}) \leq c_2 \exp(-aG^{-b} n^2) ,  \label{602}
\eea
where
\beas
& & c_1= \exp \Big( \inf\limits_n \big(ab \beta G^{1-b} S^{(b-1)\beta}n^{2-2/b} - |\tilde
R_n| \big) \Big) \  ,  \\
& & c_2= \exp \Big( \sup\limits_n \big(ab \beta G^{1-b} S^{(b-1)\beta}n^{2-2/b} + |\tilde
R_n| \big) \Big).
\eeas

Let now $n \in \bbN$ be given and let $m \in \bbN$ be such that
$m_n < m \leq m_{n+1}$, and consider $r_{m_n}^m$. We note by Lemma \ref{lem4.2}
and \eqref{258} that
 \beas
& & \log\big( r_{m_n}^m \big) = m \log \big( 1 - Gn^{-2/b} + O(n^{-4 /b}) \big) .
% \label{605}
\eeas
Since $m \leq m_{n+1} \leq C n^{\alpha } = Cn^{2 +2 /b}$ and
$ 2 +2/b  -4 /b\leq 0 $, we find using the Taylor expansion \eqref{233b} that
\beas
m \log \big( 1 - Gn^{-2/b} + O(n^{-4 /b}) \big)
= m \log  \big( 1 - Gn^{-2/b}  \big) + O(1)
\eeas
which implies
\bea
C_1 \big( 1 - Gn^{-2/b} \big)^m \leq r_{m_n}^m
\leq C_2 \big( 1 - Gn^{-2/b} \big)^m   \label{606}
\eea
for some constants $0< C_1 < C_2$.
Combining \eqref{227}, \eqref{602}, and \eqref{606} yields \eqref{230}
of Theorem \ref{th4.5}.

However, if $1 < b \leq 2 $,  in \eqref{241x}, the exponent of
$m_n^{(b-1)\beta} = n^{(b-1)\alpha \beta} = n^{2-2/b}$ is positive, although $|\tilde R_n|$ is
bounded. Instead of \eqref{602} we use
\beas
c_1 \exp(-aG^{-b} n^2 + \beta (ab)^{2\beta} n ) \leq v(r_{m_n}) \leq c_2
\exp(-aG^{-b} n^2 + \beta (ab)^{2\beta}  n ) ,
\eeas
since $abG^{1-b} =  (ab)^{2\beta}$. %%by \eqref{258}.
Lemma \ref{lem4.2} yields  for $r_{m_n}^m$
\beas
& & \log\big( r_{m_n}^m \big) = m \log \big( 1 - Gn^{-2/b} +
\beta G^2 n^{-4/b} + O(n^{-6 /b}) \big)
\roweq
m \Big( - Gn^{-2/b} + \big( \beta - \frac12 \big) G^2 n^{-4/b} + O(n^{-6 /b}) \Big)
\\
\Rightarrow & &  \log\big( r_{m_n}^m \big) =
m \log \big( 1 - Gn^{-2/b} +
\beta G^2 n^{-4/b} \big)  + O(1) ,
\eeas
since here
$mn^{-6 /b} \leq m_{n+1} n^{-6 /b} \leq  Cn^{2 +2 /b -6/b} \leq C'$.
Hence, we have
\beas
& & C_1 \big( 1 - Gn^{-2/b}  + \beta G^2 n^{-4/b} \big)^m  \leq
r_{m_n}^m
\leq  C_2 \big( 1 - Gn^{-2/b}  + \beta
G^2 n^{-4/b} \big)^m ,
\eeas
and thus \eqref{231} follows. \ \ $\Box$

\BER
\label{rem4}
It seems that the calculation of the numbers $m_n$ for  the weights
\beas
v(r)= \exp(- a/(1-r)^b)
\eeas
with our method becomes increasingly difficult for  large $b$. Technical problems
are caused  by the fact that using  an asymptotic expansion like \eqref{235} to
evaluate $v(r_m)$, more terms are required, depending on how large $b$ is.
Also, it seems that one would need a more complicated ansatz
$$
m_n= S_\alpha n^\alpha + S_{\alpha-1} n^{\alpha -1} + \ldots \ \ .
$$
The technical difficulties become obvious.
\ENR

\section{The space of multipliers $( H^\infty_v(\bbD), \ell^p ) $.}\label{secmult}

Let $A$ and $B$ be vector spaces of complex sequences containing the space of all the sequences with finitely many non-zero coordinates. The set of multipliers from $A$ into $B$ is

$$(A,B):= \{ c=(c_n) \, : \, (c_na_n) \in B \ \forall (a_n) \in A \}. $$

Given a strictly increasing, unbounded sequence $J= (m_n)_{n=0}^\infty
\subset \bbN$  and  $1 \leq p,q \leq \infty$ we denote as in \cite{BlZ}, Definition 2,
\begin{equation}
\ell^J (p,q) := \Big\{ (a_m)_{m=0}^\infty \, : \, \Big( \sum_{m=m_n+1}^{m_{n+1}}
|a_m|^p \Big)^{1/p} \in \ell_q \Big\},  \nonumber
\end{equation}
with the obvious changes when $p$ or $q$ is $\infty$. The space $\ell^J (p,q)$
is a Banach space when endowed with  the canonically defined norm. Observe
that $\ell^J(p,p) = \ell_p$. We recall the following result from \cite[Theorem 23]{BlZ} (see also \cite[Lemma 5.1]{BT}).

\begin{lemma}
\label{lem5.1}
For $1 \leq p \leq \infty$ we have
\begin{equation}
\big( \ell^J(2,\infty) , \ell^p \big) =
\ell^J(r,s)  \nonumber
\end{equation}
where (a) $r=2p/(2-p)$, $s=p$, if $1 \leq p < 2$, (b) $r=\infty$,
$s=p$, if $2 \leq p < \infty$, and (c) $r= s=\infty$, if $p = \infty$.
\end{lemma}

\BEP \label{propmult}
Let $v(r) = \exp(-1/(1-r))$, $r \in [0,1[$ and $1 \leq p \leq \infty$.
Then, the space of multipliers $\big( H_v^\infty(\bbD), \ell_p \big)$ is the set of sequences
$(\lambda_m)_{m=0}^\infty$ such that
\begin{equation}
\Big( \sum_{n=1}^\infty  \Big( \sum_{m=n^4 +1}^{(n+1)^4}
\big( |\lambda_m| e^{n^2} \Big(1  - \frac1{n^2} \Big)^{-m}  \big)^{\frac{2p}{2-p}}
\Big)^{\frac{2-p}{2}} \Big)^{\frac{1}{p}} < \infty  ,  \nonumber
\end{equation}
if $1 \leq p < 2$,
\begin{equation}
\Big( \sum_{n=1}^\infty  \Big( \max\limits_{n^4 < m \leq (n+1)^4}
 |\lambda_m| e^{n^2} \Big(1  - \frac1{n^2} \Big)^{-m}  \Big)^{p}
 \Big)^{\frac{1}{p}} < \infty  ,  \nonumber
\end{equation}
if $2 \leq p < \infty$, and
\begin{equation}
\sup\limits_{n \in \bbN } \Big( \max\limits_{n^4 < m \leq (n+1)^4}
 |\lambda_m| e^{n^2}\Big(1  - \frac1{n^2} \Big)^{-m} \Big) < \infty  ,  \nonumber
\end{equation}
if $p =  \infty$,
\ENP

Proof.
Since $\ell^p$ is a solid space, we have (cf. \cite{AS})
\begin{equation}
\big( H_v^\infty( \bbD) , \ell_p \big)=
\big( S( H_v^\infty( \bbD)) , \ell_p \big).   \nonumber
\end{equation}
Now, by Theorem  \ref{th4.5} it is easy to see that
$(\lambda_m)_{m=0}^\infty \in \big(  S( H_v^\infty( \bbD)) , \ell_p \big)$, if and
only if
\begin{equation}
\Big( \big( e^{n^2}\Big(1  - \frac1{n^2} \Big)^{-m} \Big) |\lambda_m| \big)_{m=n^4+1}^{(n+1)^4}
\Big)_{n=0}^\infty \in \big( \ell^J(2, \infty) , \ell_p\big) .  \nonumber
\end{equation}
The conclusion now follows from Lemma \ref{lem5.1}.
$\Box$

It is clear that Proposition \ref{propmult} can be extended to more weights, but we prefer to present here only this more precise formulation as an example.

\vspace{.5cm}

\noindent \textbf{Acknowledgements.} The research of Bonet was partially
supported by the projects MTM2013-43540-P and MTM2016-76647-P. This paper was completed  during the
Bonet's stay at the Katholische Universit\"{a}t Eichst\"att-Ingol\-stadt (Germany). The support of the
Alexander von Humboldt Foundation is greatly appreciated. The research of Taskinen was
partially supported by the V\"ais\"al\"a Foundation of the Finnish Academy
of Sciences and Letters.

\noindent \textbf{Authors' addresses:}%
\vspace{\baselineskip}%

Jos\'e Bonet: Instituto Universitario de Matem\'{a}tica Pura y Aplicada IUMPA,
Universitat Polit\`{e}cnica de Val\`{e}ncia,  E-46071 Valencia, Spain

email: jbonet@mat.upv.es \\

Jari Taskinen: Department of Mathematics and Statistics, University of Helsinki,
P.O.Box 68, 00014 Helsinki, Finland.

email: jari.taskinen@helsinki.fi

\end{document}